\documentclass[11pt]{article}
\usepackage{amssymb}
\usepackage{amsfonts}
\usepackage{amsmath}
\usepackage{amsopn}
\usepackage{enumerate}

\def\disp{\displaystyle}

\def\crr{\cr\noalign{\vskip2mm}}

\def\dref#1{(\ref{#1})}

\def\sqr#1#2{{\vcenter{\vbox{\hrule height.#2pt
              \hbox{\vrule width.#2pt height#1pt \kern#1pt \vrule width.#2pt}
              \hrule height.#2pt}}}}
\def\signed #1{{\unskip\nobreak\hfil\penalty50
              \hskip2em\hbox{}\nobreak\hfil#1
              \parfillskip=0pt \finalhyphendemerits=0 \par}}
\def\endpf{\signed {$\sqr69$}}
\def\3n{\negthinspace \negthinspace \negthinspace }
\def\2n{\negthinspace \negthinspace }
\def\1n{\negthinspace }

\def\={\buildrel \triangle \over =}

\def\ds{\displaystyle}

\def\ns{\noalign{\ss}}

\def\a{\alpha}
\def\b{\beta}

\def\ve{\varepsilon}

\def\f{\varphi}

\def\D{\Delta}

\def\O{\Omega}

\def\cB{{\cal B}}

\def\cF{{\cal F}}
\def\cG{{\cal G}}

\def\cK{{\cal K}}

\def\cT{{\cal T}}

\def\cW{{\cal W}}

\def\ss{\smallskip}

\def\bs{\bigskip}
\def\q{\quad}

\def\lan{\mathop{\langle}}
\def\ran{\mathop{\rangle}}

\def\cd{\cdot}

\def\deq{\mathop{\buildrel\D\over=}}

\def\({\Big (}
\def\){\Big )}
\def\[{\Big[}
\def\]{\Big]}
\def\bde{\begin{definition}}
\def\ede{\end{definition}}
\def\be{\begin{equation}}
\def\bel{\begin{equation}\label}
\def\ee{\end{equation}}
\def\bt{\begin{theorem}}
\def\et{\end{theorem}}
\def\bc{\begin{corollary}}
\def\ec{\end{corollary}}
\def\bl{\begin{lemma}}
\def\el{\end{lemma}}
\def\bp{\begin{proposition}}
\def\ep{\end{proposition}}
\def\bas{\begin{assumption}}
\def\eas{\end{assumption}}
\def\br{\begin{remark}}
\def\er{\end{remark}}
\def\ba{\begin{array}}
\def\ea{\end{array}}
\def\ed{\end{document}}

\def\square#1{\vbox{\hrule\hbox{\vrule height#1%
     \kern#1\vrule}\hrule}}
\def\rectangle#1#2{\vbox{\hrule\hbox{\vrule height#1%
     \kern#2\vrule}\hrule}}

\font\tenbb=msbm10 \font\sevenbb=msbm7 \font\fivebb=msbm5

\newfam\bbfam
\scriptscriptfont\bbfam=\fivebb \textfont\bbfam=\tenbb
\scriptfont\bbfam=\sevenbb

\def\dref#1{(\ref{#1})}
\def\disp{\displaystyle}
\def\crr{\cr\noalign{\vskip2mm}}

\DeclareMathOperator{\argmax}{argmax}
\newtheorem{lemma}{Lemma}[section]
\newtheorem{theorem}[lemma]{Theorem}
\newtheorem{remark}[lemma]{Remark}

\newtheorem{corollary}[lemma]{Corollary}

\newtheorem{definition}[lemma]{Definition}
\newtheorem{proposition}[lemma]{Proposition}
\newtheorem{assumption}[lemma]{Assumption}

\makeatletter
   
   \@addtoreset{equation}{section}
\makeatother

\begin{document}
\title{Optimal Actuator Location of  Minimum Norm  Controls
for Heat Equation with General Controlled  Domain\footnote{\small
This work was carried out with the support of the National Natural
Science Foundation of China,   and the National Research Foundation
of South Africa.}}

\author{Bao-Zhu Guo$^{a,b}$, Yashan
Xu$^{c}$, Dong-Hui Yang$^{d}$
\\
$^a${\it Academy of Mathematics and Systems Science, Academia
Sinica}\\ {\it Beijing 100190,   P.R.China}\\
$^b${\it  School  of Computer Science and Applied Mathematics}\\
 {\it University of the Witwatersrand, Wits 2050, Johannesburg, South
 Africa}\\
$^c${\it School of Mathematical Sciences,  KLMNS, Fudan University}\\
{\it Shanghai 200433, P.R.China} \\
$^d${\it School of Mathematics and Statistics, Central South University}\\
{\it Changsha 410075, P.R.China}
 }

\date{February 2, 2016}

\maketitle
\begin{center}
\begin{abstract}
In this paper, we study optimal actuator  location of the minimum
norm controls  for a multi-dimensional heat equation with control
defined in the space $L^2(\Omega\times(0,T))$.
The actuator  domain
 is time-varying  in the sense that it is only required  to
have a  prescribed Lebesgue measure for any moment. We select an
optimal actuator location so that the optimal control takes its
minimal norm over all possible actuator domains.  We build a
framework of finding the Nash equilibrium so that we can develop  a
sufficient and necessary condition to characterize the optimal
relaxed solutions for both actuator location and corresponding
optimal control of the open-loop system. The existence and
uniqueness of the optimal classical solutions are therefore
concluded. As a result, we synthesize both optimal actuator location
and corresponding optimal control into the state feedbacks which
make the  optimal solutions independent of initial data.

 \vspace{0.2cm}

 \vspace{0.2cm}


\end{abstract}

\end{center}

\section{Introduction and main results}

Different to lumped parameter systems, the location of   actuator
where optimal control optimizes  the performance in systems governed
by partial differential equations (PDEs) can often be chosen
(\cite{Morris}).  Using a simple duct model, it is shown in
\cite{Morris1} that the noise reduction performance depends strongly
on actuator location. An approximation scheme is developed in
\cite{Morris} to find  optimal location of the optimal controls for
abstract infinite-dimensional systems to minimize  cost functional
with the worst choice of initial condition. In fact, the actuator
location problem has been attracted widely  by many researchers in
different contexts but  most of them are for one-dimensional PDEs,
as previously studied elsewhere \cite{A2,D2,He,He1,p2,p3,p4}, to name
just a few. Numerical research is one of   most important
perspectives  \cite{A2,Mu1,Mu4,Mu5,Tiba}, among many others.

However, there are  few results available in  literature for
multi-dimensional PDEs. In \cite{Mu2}, a problem of optimizing the
shape and position of the damping set for internal stabilization of
a linear wave equation in $\mathbb{R}^N, N=1, 2$ is considered. The paper
\cite{Mu3} considers   a numerical approximation of null controls of
the minimal $L^\infty$-norm for a linear heat equation with a
bounded potential. An interesting study is presented in \cite{p2}
where the problem of determining a measurable subset  of
 maximizing the $L^2$ norm of the restriction of the corresponding
solution to a homogeneous wave equation on a bounded open connected
subset over a finite time interval is  addressed.  In \cite{gy2},
the shape optimal design problems related to norm optimal and time
optimal of null controlled heat equation have  been considered.
However, the controlled domains in  \cite{gy2}  are limited to some
special class of open subsets  measured by the Hausdorff metric. The
same limitations can also be found in shape optimization problems
discussed in \cite{gy0,gy}. Very recently, some optimal shape and
location problems of sensors for parabolic equations with random
initial data have  been considered in \cite{p6}.

In this paper, we consider optimal actuator location of the minimal
norm controls for a multi-dimensional internal null controllable
heat equation over an open bounded domain  $\Omega$ in $\mathbb{R}^d$
space and the duration $[0,T]$.
Our internal actuator  domains are quite general: they are varying
over all possible measurable subsets $\omega(t)$ of $\Omega$ where
$\omega(t)$  is only required to have  a prescribed measure for any
decision making moment. This work is different from \cite{p6} yet
one result (Theorem \ref{theoremM}) can be considered as a refined
multi-dimensional generalization of paper \cite{Mu5} where
one-dimensional problem is considered,  as well as a solution to a
similar but  open problem for parabolic equation mentioned in paper
\cite{p3}.

Let us first state our problem. Suppose that  $\Omega\subset\mathbb
R^d$ ($d\geq 1$) is  a non-empty bounded domain with $C^2$-boundary
$\partial \Omega$.  Let $T>0$, $\a\in(0,1)$,  and let $m(\cd)$ be
the Lebesgue measure on $\mathbb {R}^d$. Denote by \be
\cW=\left\{\omega\subset \Omega\bigm| \omega {\mbox{ is measurable
with }}m(\omega)=\alpha \cdot m(\Omega)\right\}, \ee and
\begin{equation}\label{Taiwan1}
\cW_{s,T}=\left\{w\in
L^\infty\left(\Omega\times(s,T);\{\,0,\,1\,\}\right)\bigm| m(\{x|\;
w(x,t)=1\})\equiv\a\cdot m(\Omega){\mbox{ a.e. }}t\in(s,T)\right\}.
\end{equation}
It is assumed that $a(x,t)$ is analytic in $ \Omega\times(0,T)$. For
any $s\in[0,T)$ and  $\xi\in L^2(\Omega)$,
 consider the following
controlled heat equation
\begin{equation}\label{state}
 \left\{\ba{ll}
  y_t(x, t)-\Delta y(x,t)+a(x,t)y(x,t)=w(x,t) u(x,t)~~&{\mbox{in }}\Omega\times(s,T),\\
  y(x,t)=0 &{\mbox{on }}\partial\Omega\times(s,T),\\
  y(x,s)=\xi(x)&{\mbox{in }}\Omega,
  \ea\right.
\end{equation}
where $ w\in \cW_{s,T}$ is said to be, by abuse of notation, the
actuator location,  and $u\in L^2(\Omega \times(s,T))$ is the
control. It is well known that Equation (\ref{state}) admits a
unique mild solution which is denoted by $y(\cd; w, u;s,\xi)$ or
$y((x,t); w, u;s,\xi)$ when it is necessary.

 The  minimal norm control problem can be stated as follows. For any
given  $s\in[0,T)$, $\xi\in L^2(\Omega)$, and $ w\in \cW_{s,T}$, find a minimal norm control to
solve the following optimal control problem:
$$
\mbox{ \bf Problem (NP)}_{w}^{s,\xi}:\;\;   N(
w;s,\xi)\triangleq\inf\big\{\|u\|_{L^2(\Omega\times(s,T))}\,\big|\,y((x,T);
w,u;s,\xi)=0\big\}.
$$
We want to find an optimal actuator location determined by state and
design the corresponding optimal feedback control independent of
initial data $(s,\xi)\in[0,T)\times L^2(\Omega)$. {{More precisely,
we want to find two maps:}} $\cF:[0,T)\times L^2(\Omega)\mapsto \cW$
and $\cG: [0,T)\times L^2(\Omega)\mapsto L^2(\Omega)$  so that for
any $s\in[0,T)$ and $\xi\in L^2(\Omega)$,
\be\label{1-1-4}
 \left\{\ba{ll}
  y_t(x,t)-\Delta y(x,t)+a(x,t)y(x,t)=\cF(t, y(\cdot,t)) \cG(t,y(\cdot,t))~~&{\mbox{in }}\Omega\times(s,T),\\
 y(x,t)=0 &{\mbox{on }}\partial\Omega\times(s,T),\\
  y(x,s)=\xi(x)&{\mbox{in }}\Omega,
  \ea\right.
\end{equation}
 admits a unique  mild solution $y^{\cF,\cG}(\cd;s, \xi)$
satisfying $y^{\cF,\cG}((x,T);s, \xi)=0$ and
 \be\label{1-1-5} \|u^{\cF,\cG}(s,\xi)\|_{L^2(\Omega\times(s,T))}=N(
w^{\cF,\cG}(s,\xi);s,\xi)=\inf\limits_{ w\in\cW_{s,T}}N( w; s, \xi),
\ee where \be\label{1-1-6}
 w^{\cF,\cG}(s,\xi)\deq\cF(\cdot,y^{\cF,\cG}(\cd;s,\xi))\in\cW_{s,T},
\ee
and
\be\label{1-1-7}
u^{\cF,\cG}(s,\xi)\deq\cG(\cdot,y^{\cF,\cG}(\cd;s,\xi))\in L^2(\Omega\times(s,T)).
\ee

To solve this problem, we  need to discuss the following open-loop
problem  with  $s\in[0,T)$ and $\xi\in L^2(\Omega)$ being   fixed.
In particular, we need the existence and uniqueness for optimal
classical solutions to open-loop problem. A classical optimal
actuator location of the minimal norm control problem  with respect
to $(s,\xi)$ is to
 seek  $ w^{s,\xi}\in\cW_{s,T}$  to minimize $N( w;s,\xi)$:
$$
\mbox{ \bf Problem (CP)}^{s,\xi}: \;\; \ba{l}
     \bar N(s,\xi)\deq \inf\limits_{ w\in \cW_{s,T}}  N( w;s,\xi)=N( w^{s,\xi}; s,\xi)\\
    =\inf\limits_{ w\in\cW_{s,T}}\inf\limits_{u\in L^2(\Omega\times(s,T))}\big\{\|u\|_{L^2(\Omega\times(s,T))}\,\big|\,y((x,T);w,u;\xi)=0\big\}.
    \ea
$$
If such  $ w^{s,\xi}$ exists,  we say that  $ w^{s,\xi}$ is an
optimal actuator location of the optimal minimal norm  controls with
respect to $(s,\xi)$. For  Problem (NP)$^{s,\xi}_ w$, we will apply
the duality theory in the sense of Fenchel (see, e.g.,
\cite{Fenchel,Rockafellar,Lions}), namely, we  will solve the
following
 dual problem of $(NP)^{s,\xi}_ w$:
$$
\mbox{ \bf Problem (DP)}^{s,\xi}_w: \;\;  V( w;s,\xi)\deq
\inf\limits_{z\in L^2(\Omega)}J^{s,\xi}(z; w)\deq \ds\frac{1}{2}\|
w\f(\cd;z)\|^2_{L^2(\Omega\times(s,T))}+\lan \xi,\f(s;z)\ran,
$$
 where $\f(\cd, z)$ is the solution to the following equation
 \be\label{stateRA}
 \left\{\ba{ll}
  \varphi_t(x, t)+\Delta \varphi(x,t)-a(x,t)\varphi(x,t)=0 &{\mbox{in }}\Omega\times(s,T),\\
  \varphi(x,t)=0 &{\mbox{on }}\partial\Omega\times(s,T),\\
  \varphi(x,T)=z(x)&{\mbox{in }}\Omega.
  \ea\right.
\ee Furthermore, it is derived (see Lemma \ref{lemma2.4} later) that
\be V( w;s, \xi)=-\ds\frac{1}{2}N( w; s,\xi)^2, \; \forall\;
(s,\xi)\in [0,T)\times L^2(\Omega),\, w\in\cW_{s,T}. \ee Thus \be
\ba{rl}
     &\ds\frac{1}{2}\bar N(s,\xi)^2
    = \inf\limits_{ w\in \cW_{s,T}} \ds\frac{1}{2} N( w;s,\xi)^2\\
    =&\inf\limits_{ w\in \cW_{s,T}} \left[-\inf\limits_{z\in L^2(\Omega)}\left(\ds\frac{1}{2}\| w\f(\cd;z)\|^2_{L^2(\Omega\times(s,T))}+\langle \xi,\f(s;z)\ran\right)\right]\\
    =&-\sup\limits_{ w\in \cW_{s,T}}\inf\limits_{z\in L^2(\Omega)}\left[\ds\frac{1}{2}\| w\f(\cd;z)\|^2_{L^2(\Omega\times(s,T))}+\langle \xi,\f(s;z)\ran\right].
    \ea
\ee
Therefore Problem (CP)$^{s,\xi}$ can be transformed into the following Stackelberg Problem
$$
\mbox{\bf Problem (SP)}^{s,\xi}: \;\;
 \sup\limits_{
w\in\cW_{s,T}}\inf\limits_{z\in L^2(\Omega)}\left[\ds\frac{1}{2}\|
w\f(\cd;z)\|^2_{L^2(\Omega\times(s,T))}+\langle\xi,\,\f(s;z)\rangle\right].
$$
Since $\cW_{s,T}$ lacks of compactness, it is nature to extend the
feasible set $\cW_{s,T}$ to a relaxed set $\cB_{s,T}$ (see
\dref{Guore})    to ensure the existence. But it is well
known that usually \be\label{1.1.7}
\sup\limits_{w\in\cW_{s,T}}\inf\limits_{z\in
L^2(\Omega)}\neq\sup\limits_{\beta\in \cB_{s,T}}\inf\limits_{z\in
L^2(\Omega)} \ee in the framework of game theory.

 One novelty of present work  is that the
 results derived from the  relaxed case can be returned back to the classical case.
It is difficult to verify directly if (\ref{1.1.7}) is true or not.
Our way is to prove  that any relaxed solution
 is also classical by using  a  sufficient and necessary
 condition for relaxed solutions.
 As for two-level optimization Problem (SP)$^{s,\xi}$, it is still not easy to obtain a sufficient and necessary condition.
  It is especially critical that Problem (DP)$^{s,\xi}_w$ may have  no solution in $L^2(\Omega)$
 though  Problem (NP)$^{s,\xi}_w$ always admits its solution, which
is another difficulty.  We observe keenly that in these  cases, the
Stackelberg game problem can be transformed into a Nash equilibrium
problem in a zero-sum game framework, for which a sufficient and
necessary condition  for the optimal solutions (actuator location
and control)  can be derived.

  Define
 \be\label{1-1-12}
 Z=\left\{z\in H^{-1/2}(\Omega)\bigm|\f(\cd;z)\in L^2(\Omega\times(0,T))\right\},
 \ee
 where $\f(\cd,z)$ is the solution to Equation (\ref{stateRA}) with $s=0$.
One of the main results of this  paper is     Theorem
\ref{theoremM}.
 \bt\label{theoremM}
 Let $T>0$, $\a\in(0,1)$,  and let  $a(x,t)$ be  analytic in $
\Omega\times(0,T)$. Problem (CP)$^{s,\xi}$ admits a unique solution for any $(s,\xi)\in[0,T)\times L^2(\Omega)\setminus\{0\}$.
 In
addition,
  $\bar w $ is a solution to Problem (CP)$^{s,\xi}$ if and only
if there is $ \bar z\in Z$ such that $(\bar w, \bar z)$ is a Nash
equilibrium of the following two-person zero-sum game problem: Find
$(\bar w,\,\bar z)\in\cW_{s,T}\times Z$ such that
\begin{equation}\label{gxy1-9}
\begin{array}{c}
\ns\ds\frac{1}{2}\left\|\bar w\f(\cd;\bar z)\right\|^2_{L^2(\Omega\times(s,T))}+\lan
\xi,\f(s;\bar z)\ran=\sup\limits_{w\in\cW_{s,T}}
\left[\ds\frac{1}{2}\left\| w\f(\cd;\bar z)\right\|^2_{L^2(\Omega\times(s,T))}+\lan
\xi,\f(s;\bar z)\ran\right],\\
\ns\ds\frac{1}{2}\left\|\bar w\f(\cd;\bar z)\right\|^2_{L^2(\Omega\times(s,T))}+\lan
\xi,\f(s;\bar z)\ran=\inf\limits_{z\in Z}
\left[\ds\frac{1}{2}\left\|\bar w\f(\cd; z)\right\|^2_{L^2(\Omega\times(s,T))}+\lan
\xi,\f(s;z)\ran\right].
\end{array}
\end{equation}
 \et
Another main result of this  paper is    Theorem \ref{theoremM12}.
 \bt\label{theoremM12}
 Let $T>0$, $\a\in(0,1)$,  and let  $a(x,t)$ be  analytic in $
\Omega\times(0,T)$. There are
 two maps: $\cF:[0,T)\times L^2(\Omega)\mapsto \cW$ and
$\cG: [0,T)\times L^2(\Omega)\mapsto L^2(\Omega)$  so that for any
$s\in[0,T)$ and $\xi\in L^2(\Omega)\setminus\{0\}$, Equation (\ref{1-1-4}) admits
a unique mild solution $y^{\cF,\cG}(\cd;s, \xi)$ satisfying
$y^{\cF,\cG}((x,T);s, \xi)=0$ and (\ref{1-1-5})-(\ref{1-1-7}). \et

We proceed as follows. Define
 \be\label{Guore}
 \cB_{s,T}=\left\{\beta\in L^\infty(\Omega\times(s,T);[0,1])\Bigm| \ds\int_\Omega\beta^2(x,t)\mathrm{d}x\equiv\a \cdot{ m}(\Omega) {\mbox { a.e. }}t\in(s,T)\right\}
\ee as a  relaxed set of $\cW_{s,T}$. In section 2, we discuss the
minimum norm control Problem (NP)$^{s,\xi}_\beta$ in the relaxed
case by replacing $w\in \cW_{s,T}$ with  $\beta\in \cB_{s,T}$,  and
present Problem (SP)$^{s,\xi}$ in the relaxed case. Section 3 will be
devoted to discussing   properties of solutions to Problem
(SP)$^{s,\xi}$ in the relaxed case. Finally, we prove Theorem
\ref{theoremM},  and Theorem \ref{theoremM12} is proved  by the
synthetic method.

\section{Relaxed minimum norm control problem (NP)$^{s,\xi}_\beta$}
Let $(s,\xi)\in[0,T)\times L^2(\Omega)$ be fixed. For any  $\beta\in\cB_{s,T}$.
 Consider the following system:
\be\label{stateR}
 \left\{\ba{ll}
  y_t(x, t)-\Delta y(x,t)+a(x,t)y(x,t)=\beta(x,t) u(x,t)~~&{\mbox{in }}\Omega\times(s,T),\\
  y(x,t)=0 &{\mbox{on }}\partial\Omega\times(s,T),\\
  y(x,s)=\xi(x)&{\mbox{in }}\Omega,
  \ea\right.
\ee where once again the  control $u\in L^2(\Omega\times(s,T))$, and
the solution of (\ref{stateR}) is denoted   by $y(\cd;\beta, u)$.
Accordingly, Problem $(NP)^{s,\xi}_w$ is changed into a
relaxation problem of the following:
$$
\mbox{ \bf Problem (NP)}^{s,\xi}_\beta: \;\;
 N(\beta;
s,\xi)\triangleq\inf\big\{\|u\|_{L^2(\Omega\times(s,T))}\,\big|\,y((x,T);\beta,u)=0\big\}.
$$
Let us first show the null controllability for  controlled system
(\ref{stateR}), which is deduced by building the  ``observability
inequality'' \dref{Oinequality} for system \dref{stateRA}.
 \bl\label{lemma1} For any $\beta\in\cB_{s,T}$, there exists positive constant
$C_{\beta}$ such that the solution of  \dref{stateRA} satisfies
\begin{equation}\label{Oinequality}
  \|\varphi(s;z)\|_{L^2(\Omega)}\leq C_{\beta}\left\|\beta\varphi(\cd; z)\right\|_{L^2(\Omega\times(s,T))}, \forall\; z\in L^2(\O),
\end{equation}
where $C_{\beta}$   is independent of $z\in L^2(\O)$.
 \el

\noindent {\bf Proof.} It is well known that system (\ref{stateR})
is null controllable if and only if the ``observability inequality''
\dref{Oinequality} holds for the dual  system (\ref{stateRA}). Let
$w\in\cW_{s,T}$. An  observability inequality on the  measurable set $\omega$:
\begin{equation}\label{2.7}
  \|\varphi(s;z)\|_{L^2(\Omega)}\leq \hat C_{w}\left\|w\varphi(\cd; z)\right\|_{L^2(\Omega\times(s,T))}, \forall\; z\in L^2(\O),
\end{equation}
     has been derived in \cite{AEWZ} for some $\hat C_{w}>0$.
    Now for any
$\beta\in\cB_{s,T}$, let

$$
 E=\left\{(x,t)\in\Omega\times(s,T)\bigm|\beta(x,t)\ge\sqrt{\alpha/2}\right\},
 \quad
 \lambda=\displaystyle\frac{m\left(E\right)}{m(\Omega\times(s,T))}.
$$
By
$$
\begin{array}{lll}
&~~~1\cdot m\left(\{\beta\ge \disp
\sqrt{\alpha/2}\}\right)+\alpha/2\cdot
m\left(\{\beta<\sqrt{\alpha/2}\}\right)\crr & \ge\disp
 \iint_{\{\beta\ge\sqrt{\alpha/2}\}} \beta^2(x,t)\mathrm
dx\mathrm dt+\iint_{\{\beta<\sqrt{\alpha/2}\}} \beta^2(x,t)\mathrm dx\mathrm dt\crr &=\disp
\iint_{\Omega\times(s,T)} \beta^2(x,t)\mathrm dx\mathrm dt=\alpha (T-s) m(\Omega),
\end{array}$$
here and in what follows, we denote  $\{\beta\geq \sqrt{\alpha/2}\}$ by
$\{(x,t)\in\Omega\times(s,T)\bigm| \beta(x,t)\geq
\sqrt{\alpha/2}\}$. It follows that
$$
\lambda+(1-\lambda)\alpha/2\ge \alpha .
$$
Consequently,
$
\lambda\ge\displaystyle\frac{\alpha}{2-\alpha}.
$
This  means that $E$ is not a zero-measure set.
 It then follows from (\ref{2.7}) with
$w=\chi_E$  and $\beta\ge \sqrt{\alpha/2}\,\chi_E$  that
$$
\begin{array}{rl}
\|\varphi(s;z)\|_{L^2(\Omega)}
 \leq& \hat C_{w}\left\|w  \varphi(\cdot;z)\right\|_{L^2(\Omega\times(s,T))}\le  \ds\frac{\sqrt{2}\hat C_{w}}{\sqrt{\a}}\left\| \beta\varphi(\cdot;z)\right\|_{L^2(\Omega\times(s,T))}.
\end{array}
$$
This is (\ref{Oinequality}) by taking $C_\beta=\sqrt{2}\hat C_{w}/\sqrt{\a}$.
\endpf

\subsection{Relaxed dual problem (DP)$^{s,\xi}_\beta$}
Now we present the relaxed dual problem
$$
\mbox{\bf Problem (DP)}^{s,\xi}_\beta:\;\; V(\beta;s,\xi)\deq
\inf\limits_{z\in L^2(\Omega)}J^{s,\xi}(z;\beta)\deq
\ds\frac{1}{2}\|\beta\f(\cd;z)\|^2_{L^2(\Omega\times(s,T))}+\lan
\xi,\f(s,z)\ran.
$$
Since there may have no solution in $L^2(\Omega)$ for Problem
$(DP)^{s,\xi}_\beta$,
 we need
to introduce a class of spaces $\{\bar Y_\beta,\,\beta\in\cB_{s,T}\}$. Let
\begin{equation}\label{Taiwan3}
  Y=\{\varphi(\cdot; z)|\ z\in L^2(\Omega)\}\subset L^2(\Omega\times(s,T)),
  \end{equation}
where $\varphi(\cdot; z)$ is the solution of (\ref{stateRA}) with
the initial data $z\in L^2(\Omega)$. Obviously, $Y$ is a linear
space from the linearity of  PDE  (\ref{stateRA}).

\bl\label{lemmaT} Let $Y$ be defined by \dref{Taiwan3}.  For each
$\beta\in \cB_{s,T}$, define a functional  in $Y$ by
\begin{equation*}
  F_0(\varphi)=\|\beta\varphi\|_{L^2(\Omega\times(s,T))}, \;  \forall\; \varphi\in Y.
\end{equation*}
Then $(Y, F_0)$ is a linear normed  space. We denote this normed
space by $Y_{\beta}$. \el

\noindent{\bf  Proof.} It suffices to show that
$F_0(\psi)=\|\beta\psi\|_{L^2(\Omega\times(s,T))}=0$ implies
$\psi=0$. Actually, by $F_0(\psi)=0$, it follows that
\begin{equation*}
  \sqrt{\alpha/2}\|\chi_{\{\beta\ge\sqrt{\alpha/2}\}}\psi\|_{L^2(\Omega\times(s,T))}\leq \|\beta\psi\|_{L^2(\Omega\times(s,T))}=0.
\end{equation*}
By the unique continuation (see, e.g.,  \cite{AEWZ}) for  heat
equation, we arrive at $\psi(x,t)=0$.
\endpf

Denote by
\begin{equation}\label{Taiwan5}
\overline{Y}_{\beta} =\hbox{ the completion  of the space }
Y_{\beta}.
\end{equation}
It is usually hard to characterize $\overline{Y}_{\beta}$.
However, we have the following description for $\overline{Y}_{\beta}$.
\bl\label{lemma2.2}
Let $\b\in\cB_{s,T}$,   and let $\overline{Y}_{\beta}$ be defined by \dref{Taiwan5}.  Then under an isometric
isomorphism, any element of $ \overline{Y}_{\beta}$ can be
expressed as a function $\hat\f\in C([s,T); L^2(\Omega))$ which
satisfies (in the sense of weak solution)
\begin{equation}\label{stateRAguo}
 \left\{\ba{ll}
  \hat{\varphi}_t(x, t)+\Delta \hat{\varphi}(x,t)-a(x,t)\hat{\varphi}(x,t)=0~~&{\mbox{\rm in }}\Omega\times(s,T),\\
  \hat{\varphi}(x,t)=0 &{\mbox{\rm on }}\partial\Omega\times(s,T),
  \ea\right.
\end{equation}
and $\beta\hat\f=\lim\limits_{n\rightarrow\infty}\beta\varphi(\cdot;
z_n)$ in
$L^2(\Omega\times(s,T))$ for some sequence  $\{z_n\}\subset L^2(\Omega)$, where $\varphi(\cdot; z_n)$ is the solution
of \dref{stateRA} with initial value $z=z_n$. \el
\noindent{\bf
Proof.}
 Let $\overline{\psi}\in(\overline{Y}_{\beta},\bar F_0)$, where $(\overline{Y}_{\beta},\bar F_0)$ is the completion of $(Y_{\beta},F_0)$.
 By  definition, there is a sequence
$\{z_n\} \subset L^2(\Omega)$   such that
\begin{equation*}
  \bar F_0(\varphi(\cdot; z_n)-\overline{\psi})\rightarrow 0,
\end{equation*}
from which, one has
\begin{equation*}
  F_0\left(\varphi(\cdot; z_n)-\varphi(\cdot; z_m)\right)=\bar F_0
(\varphi(\cdot; z_n)-\varphi(\cdot; z_m))\rightarrow 0 \hbox{ as } n, m \to \infty.
\end{equation*}
In other words,
\begin{equation}\label{10.6.1}
  \|\beta\varphi(\cdot; z_n)-\beta\varphi(\cdot; z_m)\|_{L^2(\Omega\times(s,T))}\rightarrow 0 \hbox{ as } n, m \to \infty.
\end{equation}
Hence, there exists $\hat\psi\in L^2(\Omega\times(s,T))$ such that
\begin{equation}\label{2.9}
\beta \varphi(\cdot; z_n)\rightarrow \hat\psi\;\;\mbox{strongly
in}\;\; L^2(\Omega\times(s,T)).
\end{equation}
Let $\{T_k\}\subset (s,T)$ be such that $T_k \nearrow T$. i.e. $T_k$ is
strictly monotone increasing  and converges to $T$. Denote $\varphi_n\equiv\varphi(\cdot;
z_n)$.

(a).  For $T_1$, by the observability inequality (\ref{Oinequality}),
and (\ref{10.6.1}),
\begin{eqnarray*}\label{2.10}
\begin{array}{rl}
\ns&\|\varphi(T_2; z_n)\|_{L^2(\Omega)} \leq C(1)\|\beta\varphi(\cdot;
z_n)\|_{L^2(T_2,T; L^2(\Omega))}\nonumber\\
\ns \leq& C(1)\|\beta\varphi(\cdot;
z_n)\|_{L^2(\Omega\times(s,T))} \leq C(1)\sup\limits_m\|\beta\varphi(\cdot;
z_m)\|_{L^2(\Omega\times(s,T))}, ~\forall\;
n\in \mathbb{N},
\end{array}\end{eqnarray*}
Hence,  there exists a subsequence $\{\varphi_{1n}\}$ of $\{\varphi_n\}$
 and $\varphi_{01}\in L^2(\Omega)$ such that
\begin{equation*}
   \varphi_{1n}(T_2)=\varphi(T_2; z_{1n})\rightarrow \varphi_{01} \mbox{ weakly in } L^2(\Omega).
\end{equation*}
This together with the fact:
\begin{equation*}
 \left\{\ba{ll}
  (\varphi_{1n})_t(x, t)+\Delta \varphi_{1n}(x,t)-a(x,t)\varphi_{1n}(x,t)=0~~&{\mbox{in }}\Omega\times(s,T_2),\\
  \varphi_{1n}(x,t)=0 &{\mbox{on }}\partial\Omega\times(s,T_2),\\
  \varphi_{1n}(x,T_2)=\varphi(T_2; z_{1n})&{\mbox{in }}\Omega,
  \ea\right.
\end{equation*}
shows that there exists $\psi_1\in L^2(s,T_2; L^2(\Omega))\cap
C([s,T_2-\delta];L^2(\Omega))$ for all $\delta>0$,  satisfies
\begin{equation*}\label{10.7.4guo1}
 \left\{\ba{ll}
  (\psi_1)_t(x, t)+\Delta \psi_1(x,t)-a(x,t)\psi_1(x,t)=0~~&{\mbox{in }}\Omega\times(s,T_2),\\
  \psi_1(x,t)=0 &{\mbox{on }}\partial\Omega\times(s,T_2),\\
  \psi_1(x,T_2)=\varphi_{01}(x)&{\mbox{in }}\Omega,
  \ea\right.
\end{equation*}
and  for  all  $\delta\in(0,T_2)$,
\begin{equation*}
  \varphi_{1n}\rightarrow \psi_1  \mbox{ strongly in } L^2([s,T_2]; L^2(\Omega))\cap C([s,T_2-\delta]; L^2(\Omega)).
\end{equation*}
In particular,
\begin{equation}\label{10.7.1}
   \varphi_{1n}\rightarrow \psi_1 \mbox{ strongly in } L^2([s,T_2]; L^2(\Omega))\cap C([s,T_1];
   L^2(\Omega)),
\end{equation}
and
\begin{equation}\label{10.7.2}
  \beta\varphi_{1n}\rightarrow \beta\psi_1  \mbox{ strongly in } L^2([s,T_2];
  L^2(\Omega)).
\end{equation}
These together with  (\ref{2.9}) and (\ref{10.7.2}) yield
\begin{equation*}\label{10.7.3}
  \beta\psi_1=\hat\psi \mbox{ in } L^2([s,T_1]; L^2(\Omega)).
\end{equation*}

(b). Along the same way as (a),  we can find a subsequence
$\{\varphi_{2n}\}$ of $\{\varphi_{1n}\}$, and $\psi_2\in
L^2([s,T_3];L^2(\Omega))\cap C([s,T_3-\delta];L^2(\Omega))$ for all
$\delta>0$  satisfying
\begin{equation*}\label{10.7.4guo2}
 \left\{\ba{ll}
  (\psi_2)_t(x, t)+\Delta \psi_2(x,t)-a(x,t)\psi_2(x,t)=0~~&{\mbox{ \rm in }}\Omega\times(s,T_3),\\
  \psi_2(x,t)=0 &{\mbox{ \rm on }}\partial\Omega\times(s,T_3),
  \ea\right.
\end{equation*}
and
\begin{equation*}
   \varphi_{2n}\rightarrow \psi_2 \mbox{ strongly in } L^2([s,T_3]; L^2(\Omega))\cap C([s,T_2]; L^2(\Omega)).
\end{equation*}
This,  together with (\ref{10.7.1}), leads to
\begin{equation*}
  \psi_2|_{[s,T_1]}=\psi_1,
\end{equation*}
and
\begin{equation*}\label{10.7.3*}
  \beta\psi_2=\hat\psi \mbox{ in } L^2([s,T_2]; L^2(\Omega)).
\end{equation*}

(c). Similarly  to  (a) and (b), we can find a  sequence
$\{\psi_k\}$ which satisfies, for each $k\in \mathbb{N}^+$, that

\begin{itemize}

\item  $\psi_k\in L^2([s,T_{k+1}]; L^2(\Omega))\cap C([s,T_k];
L^2(\Omega))$;

\item $  \psi_{k+1}|_{[s,T_k]}=\psi_k$;

\item $\psi_k$ satisfies
\begin{equation*}\label{10.7.4guo3}
 \left\{\ba{ll}
  (\psi_k)_t(x, t)+\Delta \psi_k(x,t)-a(x,t)\psi_k(x,t)=0 &{\mbox{ in }}\Omega\times(s,T_{k+1}),\\
  \psi_k(x,t)=0 &{\mbox{on }}\partial\Omega\times(s,T_{k+1}).
  \ea\right.
\end{equation*}

 \item $  \beta\psi_k=\hat\psi \mbox{ in } L^2([s,T_k]; L^2(\Omega))$.
\end{itemize}
Define
\begin{equation*}
  \psi(\cdot,t)=\psi_k(\cdot,t), \;  t\in [s,T_k].
\end{equation*}
Then,  $\psi(x,t)$ is  well defined on $[s,T)$,  which satisfies
$\psi\in L^2([s,T]; L^2(\Omega))\cap C([s,T); L^2(\Omega))$,
\begin{equation*}\label{10.7.4guo4}
 \left\{\ba{ll}
  \psi_t(x, t)+\Delta \psi(x,t)-a(x,t)\psi(x,t)=0~~&{\mbox{in }}\Omega\times(s,T),\\
  \psi(x,t)=0 &{\mbox{on }}\partial\Omega\times(s,T),
  \ea\right.
\end{equation*}
and
$$
   \beta\psi=\hat\psi=\lim\limits_{n\rightarrow\infty}\beta\varphi(\cdot; z_n).
   $$
Under   an  isometric isomorphism,   we can say
$\overline{\psi}=\psi$. This complete the proof of the lemma.
\endpf

\bs

We  define an operator $\cT: Y\rightarrow L^2(\Omega)$ by \be
\label{operator1} \cT(\varphi(\cd; z))=\varphi(s; z), \; \forall\;
z\in L^2(\Omega), \ee
which is well-defined because  $Y\subset C([0,T]; L^2(\Omega))$.
Define an operator $\cT_{\beta}: \beta\overline{Y}_{\beta}\rightarrow L^2(\Omega)$ by \be\label{operator2} \cT_{\beta}(\beta\psi)=\psi(s),  \; \forall\; \psi\in \overline{Y}_{\beta}. \ee
 By lemma \ref{lemma2.2}, the operator   $\cT_{\beta}$ is  also well-defined. In addition,  it follows from the observability  inequality claimed by  Lemma
\ref{lemma1} that  the  linear operator $\cT_{\beta}$ is bounded.

\bl \label{lemma2.3} If  $\beta\in\cB_{s,T}$, then
the operator $\cT_{\beta}$ defined by \dref{operator2}   is
compact. \el {\bf  Proof.} By  the observability inequality claimed
by Lemma \ref{lemma1}, it follows that
 the operator
${\beta\overline{Y}_{\beta}}\rightarrow L^2(\Omega)$ defined by
$$
\beta\psi(\cd, \cd)\rightarrow \psi(\cd, (T+s)/2), \; \forall\; \psi\in
\overline{Y}_{\beta}
$$
is bounded. Also by the property of heat equation, the operator
defined by
$$
\f(\cd, (T+s)/2)\rightarrow \f(\cd, s), \; \forall\;
\f\in\overline{Y}_{\beta}
$$
is compact.  As a composition operator from the above two operators,
$\cT_{\beta}$ is compact as well.
\endpf

\vspace{3mm}

Now we tune to  discuss the solution to Problem (DP)$^{s,\xi}_\beta$ with extended domain. From  the notation of $\cT_\beta$, we could rewrite the functional $J^{s,\xi}(\cdot;\beta)$ in Problem (DP)$^{s,\xi}_\beta$ as follows:
$$
J^{s,\xi}(\zeta;\beta)=\ds\frac{1}{2}\|\zeta(\cd)\|^2_{L^2(\Omega\times(s,T))}+\lan\xi,\cT_\beta(\zeta)\ran, \; \forall\; \psi\in Y.
$$
 Let us expand the domain of $J^{s,\xi}(\cdot;\beta)$ as follows:
$$
\hat J^{s,\xi}_\beta(\zeta)=\ds\frac{1}{2}\|\zeta(\cd)\|^2_{L^2(\Omega\times(s,T))}+\lan\xi,\cT_\beta(\zeta)\ran
\quad{\mbox{for any }}\zeta\in \beta\overline Y_\beta
$$
and denote
$$
\mbox{\bf Problem ($\widehat {DP}$)}^{s,\xi}_\beta:\;\; V(\beta;s,\xi)=
\inf\limits_{\zeta\in\beta\overline Y_\beta}\hat J^{s,\xi}_\beta(\zeta)\deq
\ds\frac{1}{2}\|\zeta(\cd)\|^2_{L^2(\Omega\times(s,T))}+\lan\xi,\cT_\beta(\zeta)\ran
$$
In the above, the first equation holds from Lemma \ref{lemma2.2} and the continuity of $\tau_\beta$.

\subsection{Relationship between Problems (NP)$^{s,\xi}_\beta$ and   ($\widehat {\mbox{DP}}$)$^{s,\xi}_\beta$}

In this subsection, we  present two properties on the relationship
between Problems $(NP)_{p,\beta}$ and  ($\widehat {DP}$)$^{s,\xi}_\beta$.

\begin{lemma}\label{lemma2.4}
Let $s\in(0,T)$, $\xi\in L^2(\Omega)\setminus\{0\}$, and $\beta\in\cB_{s,T}$.
Then
Problem ($\widehat {DP}$)$^{s,\xi}_\beta$ admits a unique nonzero solution in $\beta\overline Y_\beta$, denoted by $\bar\zeta$,
and the control defined by
 \be\label{optimal control}
 \bar u(x,t)= \bar\zeta(x,t),\quad
  (x,t)\in\Omega\times(s,T)  \hbox{ a.e.}
 \ee
 is an optimal control to Problem (NP)$^{s,\xi}_\beta$. Moreover,
 \begin{equation}\label{value}
V(\beta;s,\xi)=-\frac{1}{2}N(\beta;s,\xi)^2.
\end{equation}
\end{lemma}
{\bf Proof.}
 Since  $L^2(\Omega\times(s,T))$ is  reflexible.
Thus, $\beta\overline Y_{\beta}$, as a  closed subspace of $L^2(\Omega\times(s,T))$, is
also reflexible. Meanwhile, one can directly check that
$\hat J^{s,\xi}(\cd,\beta)$ is strictly convex and coercive  in $\beta \overline Y_\beta$.
Hence, $\hat J^{s,\xi}(\cd,\beta)$ has  a unique minimizer $\bar\zeta$.  It follows from
the unique continuity of heat equations, the map from $\overline Y_\beta$ to $\beta\overline Y_\beta$ is one-to-one. Thus there is unique $\bar\psi\in\overline Y_\beta$ such that
$\bar\zeta=\beta\bar\psi$.

We prove
 \be\label{2.12}
 \bar\zeta\neq 0 \hspace{5pt}{\mbox{and }}\bar\psi\neq0\hspace{5pt} \mbox{ in } L^2(\Omega\times(s,T)).
 \ee
Indeed, if this is not true, then  it must hold that $V(\beta;s,\xi)=0$.
We claim that $\{\f(s,z)| z\in L^2(\Omega)\}$ is dense in
$L^2(\Omega)$.  Once the claim holds, there is $z\in L^2(\Omega)$ such that
$\langle\xi, \f(s,z)\ran<0$ because $\xi\neq\{0\}$.
But,
$$0=V(\beta;s,\xi)\le J^{s,\xi}(\ve z,\beta)= \ds\frac{1}{2}\ve^2\|\beta\f(\cd;z)\|^2_{L^2(\Omega\times(s,T))}+\ve\lan
\xi,\f(s,z)\ran<0.
$$
where the  last inequality holds as  $\ve>0$ is small enough.

Now we show that $\{\f(s,z)| z\in L^2(\Omega)\}$ is dense in
$L^2(\Omega)$. Recalling the dual system  (\ref{stateRA}), we define
the operator $L$ in $L^2(\Omega)$ by
$$
Lz=\varphi(s,z) {\mbox{ for any }} z\in L^2(\Omega).
$$
Notice that
$$\begin{array}{c}
\{\varphi(s,z)| z\in L^2(\Omega)\} \hbox{ is dense in } L^2(\Omega)
\Leftrightarrow  \overline{{\cal R}(L)}=L^2(\Omega)
\Leftrightarrow {{\cal N}(L^*)}=\{0\},
\end{array}
$$
where the last equivalence holds because of $\overline{{\cal
R}(L)}={{\cal N}(L^*)}^\bot$. For any $\hat z\in L^2(\Omega)$,
consider the following equation:
$$
\left\{
\begin{array}{l}
\hat \varphi_t(t)-\triangle \hat \varphi(t)+a(T-t)\hat\varphi(t)=0,\\
\hat\varphi(s)=\hat z.
\end{array}
\right.
$$
A direct verification shows that
$$
L^*(\hat z)=\hat\varphi(T).
$$
By the backward uniqueness for heat equation, we have  ${\cal
N}(L^*)=\{0\}$, and this leads to  (\ref{2.12}).

  Now, we show that the control defined by (\ref{optimal control}) is  optimal to Problem (NP)$^{s,\xi}_\beta$.
  Since $\bar\zeta(x,t)$ is optimal,  we have
 \be\label{2.17}
 \langle \bar u,\zeta\rangle_{L^2(\Omega\times(s,T))} + \lan\xi,\cT_\beta(\zeta)(s)\ran=0,\; \forall\; \zeta\in\beta \overline{Y}_{\beta}.
  \ee
 Taking   $\zeta=\beta\f(\cd; z)$ for any $z\in L^2(\Omega)$ in \dref{2.17}, a straightforward calculation shows that
$$
 y(T; \beta,\bar u)=0.
$$
If  $\hat u(\cdot,\cdot)$ satisfies
 \be\label{2.14} y(T; \beta,\hat u)=0,
 \ee
  we will show that \be\label{2.15} \|\bar u\|_{L^2(\Omega\times(s,T))}\le
\|\hat u\|_{L^2(\Omega\times(s,T))}, \ee from which we see that
$\bar u(\cdot,\cdot)$ is an optimal  solution to Problem
(NP)$^{s,\xi}_\beta$.

Now,  we prove (\ref{2.15}). By (\ref{2.14}),
$$
-\langle \xi,\f(s,z)\rangle=\langle y(T; \beta,\hat
u),z\rangle-\langle \xi,\f(s; z)\rangle=\langle\f(\cdot;z),\beta\hat u(\cdot)
\rangle_{L^2(\Omega\times(s,T))}, \; \forall\; z\in L^2(\Omega),
$$
which is rewritten as
$$ -\lan\xi,\cT_\beta(\zeta)(s)\ran=\lan \hat
u,\zeta\rangle_{L^2(\Omega\times(s,T))}, \; \forall\; \zeta\in  \beta Y_{\beta}.
$$
By the density argument, the above still holds for any $\zeta\in \beta \overline{Y}_{\beta}.
$
It then follows from (\ref{2.17}) that
$$
\lan\bar u,\zeta\rangle_{L^2(\Omega\times(s,T))}= \lan \hat u,\zeta\rangle_{L^2(\Omega\times(s,T))}, \; \forall\; \xi\in \beta
\overline{Y}_{\beta}.
$$
Taking $\zeta=\bar\zeta$ in  above equality, we have \be\label{2.20}
\lan\bar u,\bar\zeta\rangle_{L^2(\Omega\times(s,T))}= \lan \hat u,\bar\zeta\rangle_{L^2(\Omega\times(s,T))}. \ee On the
other hand, it follows from (\ref{optimal control}) that
\be\label{2.21} \left\|\bar
u\right\|^2_{L^2(\Omega\times(s,T))}=\left\|\bar
\zeta\right\|^2_{L^2(\Omega\times(s,T))} = \lan\bar\zeta, \bar
u\rangle_{L^2(\Omega\times(s,T))}. \ee
By (\ref{2.21}) and (\ref{2.20}),
$$
\begin{array}{ll}
\|\bar u\|^2_{L^2(\Omega\times(s,T))} &=\disp \langle \bar u,\bar \zeta\rangle_{L^2(\Omega\times(s,T))}
= \langle \hat u, \bar \zeta\rangle_{L^2(\Omega\times(s,T))}\crr
  & \le\|\hat u\|_{L^2(\Omega\times(s,T))}\cdot \|\bar\zeta\|_{L^2(\Omega\times(s,T))}=\|\hat u\|_{L^2(\Omega\times(s,T))}\cdot\|\bar u\|_{L^2(\Omega\times(s,T))}.
\end{array}
$$
The inequality   $
 \| \bar u \|_{L^2(\Omega\times(s,T))}\le\| \hat
 u\|_{L^2(\Omega\times(s,T))} $ then follows immediately   because  $\bar u\neq 0$.

With a straightforward calculation, we can obtain
$$
V(\beta; s,\xi)=\ds\frac{1}{2}\|\bar\zeta\|^2_{L^2(\Omega\times(s,T))}+
 \lan\xi,\cT_\beta(\bar\zeta)(s)\ran.
$$
This, together with \dref{2.17}, \dref{optimal control} and the optimality of
$\bar u$, implies \dref{value}.
\endpf

Now we present  relaxed optimal actuator location of the minimal
norm control problem  with respect to $(s,\xi)$:
\begin{equation}\label{wits2}
\mbox{\bf Problem (RP)}^{s,\xi}:\;\;
\inf\limits_{\beta\in\cB_{s,T}}\inf\limits_{u\in
L^2(\Omega\times(s,T))}\big\{\|u\|_{L^2(\Omega\times(s,T))}\,\big|\,y(T;\omega,u;y_0)=0\big\}.
\end{equation}
 By the same argument  in  Section 1, Problem $(RP)^{s,\xi}$ is equivalent  to the following Problem (SP)$^{s,\xi}$ in the relaxed case.
\begin{equation}\label{wits3}
\mbox{\bf Problem (RSP)}^{s,\xi}:\;\;
\sup\limits_{\beta\in\cB_{s,T}}\inf\limits_{z\in
L^2(\Omega)}\left[\ds\frac{1}{2}\|\beta\f(\cd;z)\|^2_{L^2(\Omega\times(s,T))}+\langle\xi,\,\f(s;z)\rangle\right].
\end{equation}

\section{Relaxed Stackelberg game problem}

 Let us recall some basic facts of the two-person zero-sum
game problem. There are two players:  Emil and Frances. Emil takes
his strategy $x$ from his strategy set $E$ and Frances takes his
strategy $y$ from his strategy set $F$.  Let $f:\,E\times F$ be  the
index cost function. Emil  wants to minimize the function $F$ while
Frances wants to maximize $F$. In the framework of two-person
zero-sum game,  any solution to $\inf\limits_{x\in
E}\,\sup\limits_{y\in F}f(x,y)$ is called a conservative strategy of
Emil while any solution to $\sup\limits_{y\in F}\,\inf\limits_{x\in
E}f(x,y)$ is called a conservative strategy of Frances. For a game
problem, the  Nash equilibrium is the most
 important concept.

\bde \label{De.guo1*} Suppose that $ E$ and $F$ are  strategy sets of
Emil and Frances, respectively. Let $f:\,E\times F\mapsto
\mathbb{R}$ be an index cost functional. We call $(\bar x,\bar y)\in
 E\times   F$  a Nash equilibrium if,
$$
f(\bar x, y)\le f(\bar x,\bar y)\le f(x,\bar y), \quad\forall\;
x\in E, y\in  F.
$$
\ede

The following  result is well known, see, for instance, Proposition
8.1 of  \cite[p.121]{Aubin}.  It connects the Stackelberg
equilibrium with the Nash equilibrium.

 \bp \label{prpoperty3.18} The following
conditions are equivalent.
\begin{enumerate}[(i)]
\item $(\bar x,\bar y)$ is a Nash equilibrium;
\item $V^+=V^-$ and  $\bar x$ is a conservative strategy of Emil  ( equivalently, the following equation holds):
$$
V^+\deq\inf_{x\in   E}\sup_{y\in   F} f(x,y)=\sup\limits_{y\in   F} f(\bar x,y),
$$
 and $\bar y$ is a conservative strategy of Frances (equivalently, the following equation holds):
$$
V^-\deq\sup_{y\in   F}\inf_{x\in  E} f(x,y)= \inf\limits_{x\in   E} f(x,\bar y).
$$
\end{enumerate}
When $V^+=V^-$, we say that the game problem attains its value at $V^+$. \ep

Notice that Problem (RSP)$^{s,\xi}$ is a typical  Stackelberg game problem and  we will discuss it in the framework of two-person zero-sum game theory.
Let
\be
 \cB^2_{s,T}=\left\{b\in L^\infty(\Omega\times(s,T);[0,1])\Bigm| \ds\int_\Omega b(x,t)\mathrm{d}x\equiv\a \cdot{ m}(\Omega) {\mbox { a.e. }}t\in(s,T)\right\}
\ee and define an  index cost functional  by
 \begin{equation}\label{Taiwan8}
 \ba{r}
F(b,\psi)=-\ds\frac{1}{2}\displaystyle\iint_{\Omega\times(s,T)}
b(x,t)\psi^2(x,t)\mathrm dx\mathrm dt-\langle \xi,\,\psi(s)\rangle,
\forall\; \psi\in Y, \, b\in \cB^2_{s,T}.
 \ea\end{equation}
We assume that
 Emil who  controls the function $ b\in \cB^2_{s,T}$  wants to minimize $F$
and likewise,
Frances who
 controls the function $\psi\in Y$ wants to maximize
 $F$.
Then Problem (RSP)
 has the
following equivalent form:
\begin{equation}\label{wits4}
{\bf Problem (RSP1):} \;\;   V^+\deq\inf\limits_{
b\in\cB^2_{s,T}}\,\sup\limits_{\psi\in Y}F(b,\psi)
=\inf\limits_{ b\in\cB^2_{s,T}}\,\sup\limits_{\psi\in \overline
Y_\beta} F(b,\psi)
\end{equation} with $\beta=\ds\sqrt{
b}$.

\bt\label{theoremM1}
 Problem (RSP1)   admits a solution in
$\cB^2_{s,T}$.
\et
\noindent{\bf Proof.} For any $\psi\in Y$, it is clear that
  the functional $F(\cd,\psi)$ is linear and hence it is weakly*
lower semi-continuous. Let $X=L^\infty(\Omega)$ be  equipped with
the weak* topology. Then  $F(\cd,\psi)$ is lower semi-continuous
under the topology of $X$. If we denote
$$
\hat F(b)=\sup\limits_{\psi\in Y}F(b,\psi), \forall\;
b\in\cB^2_{s,T},
$$
then $\hat F(b)$ is also lower semi-continuous. In addition, since
   $\cB^2_{s,T}$ is compact  under the topology of $X$,
there exists at least one solution solving
$\inf\limits_{b\in\cB^2_{s,T}} \hat F(b)$. Therefore, the game
Problem (RSP1)    admits a solution in $\cB^2_{s,T}$.\endpf

\vspace{-0.3cm}

\subsection{Value attainability for   zero-sum game}
   In this subsection, we
will make use of the game theory  to discuss  value attainability
for above two-person zero-sum game. More precisely, denote by
\begin{equation}\label{wits5}
{\bf Problem (RSP2):}\;\;   V^-\deq\sup\limits_{\psi\in
Y}\,\inf\limits_{ b\in\cB^2_{s,T}}F(b,\psi).
 \end{equation}
 Once $V^+=V^-$, we say that the above two-person zero-sum game attains
its value. Furthermore, it is possible to characterize the
conservative strategy of Frances (solutions to Problem (RSP1)) by
using Proposition \ref{prpoperty3.18}.  To this end, we introduce an
intermediate value $\hat V$ and prove successively that $V^-=\hat V$
under topological assumptions,  and that $\hat V=V^+$ under
convexity assumptions.

We denote by $\cK$  all the finite subsets of $Y$. For any $K\in
\cK$, set
\begin{equation}\label{Taiwan9}
V_K=\inf\limits_{ b\in\cB^2_{s,T}}\sup\limits_{\psi\in
K}F(b,\psi), \; \; \hat V\deq\inf\limits_{K\in\cK}
V_K=\sup\limits_{K\in\cK}\inf\limits_{ b\in\cB^2_{s,T}}\sup\limits_{\psi\in
K}F(b,\psi).
\end{equation}
Then
 \be\label{3.10} V^-\le \hat V\le V^+. \ee

\bl\label{lemma1-3-3}~Let $V^+$ and $\hat V$ be defined by (RSP1) and (RSP2)
respectively. Then \be V^+ =\hat V. \ee \el

\noindent {\bf  Proof.} For any $\psi\in Y$, it is clear that the
  functional $F(\cd,\psi)$ is sequentially weakly*
lower semi-continuous.
 Furthermore,
for any $K=\{\psi_1, \psi_2,\dots, \psi_n\}\in \cK$, functional
$\sup\limits_{\psi\in K}F(\cd,\psi)$ is also sequentially weakly*
lower semi-continuous. This, together with the compactness of
$\cB^2_{s,T}$, implies that
 there is $ b_K\in\cB^2_{s,T}$ such
that
$$
\sup\limits_{\psi\in K}F( b_K , \psi)=\inf\limits_{ b\in\cB^2_{s,T}}\sup\limits_{\psi\in K}F( b, \psi).
$$
It then follows from the definition of $\hat V$ that \be\label{3.12}
F( b_K , \psi)\le \sup\limits_{\hat\psi\in K}F( b_K ,\hat
\psi)=\inf\limits_{ b\in\cB^2_{s,T}}\sup\limits_{\hat\psi\in K}F( b,
\hat\psi)\le\sup\limits_{\hat K\in\cK}\inf\limits_{
b\in\cB^2_{s,T}}\sup\limits_{\hat\psi\in \hat K}F( b,\hat\psi)= \hat
V, \quad\forall\; \psi\in K. \ee If we
 denote by
$$
S_\psi\deq\left\{ b\in\cB^2_{s,T}\bigm| F( b, \psi)\le \hat V\right\}
$$
for any $\psi\in Y$, then  it follows from (\ref{3.12}) that  $
b_K\in\cap_{\psi\in K}S_\psi$ and hence
\be\label{3.13}
 \bigcap_{\psi\in K}S_\psi\neq\emptyset {\mbox{ for any }}K\in\cK.
  \ee In addition, since    $F(\cd,\psi)$ is
weakly* lower semi-continuous,
  $S_\psi$ is  weakly* closed in $L^\infty(\Omega\times(s,T))$ as well. In other words,
 $S_\psi$ is   closed under the
weak* topology  of $L^\infty(\Omega\times(s,T))$.   We  claim that
\be \bigcap_{\psi\in Y} S_\psi\neq \emptyset. \ee Indeed, if the
above condition fails, then $ \bigcup_{\psi\in Y}
\cB^2_{s,T}\setminus S_\psi=\cB^2_{s,T}. $ It follows from the
compactness of $\cB^2_{s,T}$ that there is $\hat K\in\cK$ such that
$$
\bigcup_{\psi\in \hat K} \cB^2_{s,T}\setminus S_\psi=\cB^2_{s,T}.
$$
This contradicts to (\ref{3.13}).  Select
 $\bar b$ in the set $\bigcap_{\psi\in Y} S_\psi$.
 Then
$$
\sup_{\psi\in Y}F(\bar b, \psi)\le \hat V,
$$
and  so
\begin{equation*}
  \inf_{b\in\cB^2_{s,T}}\sup\limits_{\psi\in Y}F( b, \psi)\le \hat V.
\end{equation*}
This, together with (\ref{3.10}), completes the proof of  the
lemma.\endpf

\vspace{3mm}

The following Proposition \ref{prpoperty3.9} is  Proposition 8.3 of
\cite[p.132]{Aubin}. \bp \label{prpoperty3.9} Let $\hat E$ and $
\hat F$ be two convex sets and  let the function $f(\cdot,\cdot) $
be defined in $\hat E\times\hat F$. Let $\mathcal{F}$ be the set of
all finite subsets of $\hat F$ and
$$
    \hat V
= \sup_{K\in \cF}\inf_{x\in\hat E}\sup\limits_{\psi\in K}f(x,y),\;
    V^-
= \sup_{y\in \hat F}\inf_{x\in\hat E}f(x,y).
$$
Suppose that  a)~ for any $y\in \hat F$, $x\rightarrow f(x,y)$ is
convex;  and b)~for any $x\in \hat E$, $x\rightarrow f(x,y)$ is
concave. Then $\hat V=V^-$. \ep

 \bt \label{theorem3.12} Let $V^+$ and $\hat V$ be defined by (RSP1) and (RSP2), respectively.
  Then
 \be\label{1-3-9}
  V^- = V^+. \ee   \et
 \noindent {\bf Proof.} It is clear that  both $\cB^2_{s,T}$ and $Y$ are convex. We can verify directly that
 the functional $F( \cd,\psi)$ is linear and hence  is convex for any
and $\psi\in Y$. In addition,  the functional $F( b, \cd)$ is
concave for any $ b\in\cB^2_{s,T}$. Thus  $\hat V=V^-$ in terms of
Proposition \ref{prpoperty3.9}. The equality (\ref{1-3-9}) is then
derived by applying Lemma   \ref{lemma1-3-3}.
 This completes the proof of the lemma.
\endpf

\subsection{Nash equilibrium}\label{Se3.2.3}
The value  attainability  for a  given two-person zero-sum game is a
necessary condition to the existence of the Nash equilibrium. To
discuss further about the solution to the Stackleberg game Problem
(RSP1)    or equivalently  Problem (RSP)$^{s,\xi}$,  we need to
discuss another  Stackleberg game Problem (RSP2), in other words, we
should discuss the following problem: \be\label{gxy3.38}
 \inf\limits_{\psi\in
Y}\sup\limits_{ b\in\cB^2_{s,T}}\left[\ds\frac{1}{2}\int^T_s\int_\Omega b(x,t)\psi(x,t)^2\mathrm
dx\mathrm dt+\lan \xi,
\psi(s)\ran \right]. \ee Define a non-negative nonlinear functional
on $Y$ by
\begin{equation}\label{guo20}
N\hspace{-1mm}F(\psi)=\sup\limits_{ b\in\cB^2_{s,T}}\left(\ds\int^T_s\int_\Omega b(x,t)\psi(x,t)^2\mathrm dx\mathrm dt\right)^{1\over 2},
\; \forall\; \psi\in Y.
\end{equation}

\bl\label{Le.guo1} Let $N\hspace{-1mm}F(\cdot)$ be the functional
defined by \dref{guo20}. Then $N\hspace{-1mm}F(\cdot)$ is a norm for the
space $Y$ defined by \dref{Taiwan3}. \el

\noindent {\bf Proof.} It is clear that
$$
N\hspace{-1mm}F(\psi)\ge 0, \; \forall\; \psi\in Y \hbox{ and
}\psi=0\Rightarrow N\hspace{-1mm}F(\psi)= 0.
$$
By the relation between $\cB_{s,T}$ and $\cB_{s,T}^2$,
$$
N\hspace{-1mm}F(\psi)=\sup\limits_{\beta\in\cB_{s,T}}\|\beta\psi\|_{L^2(\Omega\times(s,T))}.
$$
Furthermore, if $N\hspace{-1mm}F(\psi)= 0$, then $ \beta\psi=0$ for
any $\beta\in \cB_{s,T}$. Take
$$
\hat\beta(x,t)\equiv\chi_{\omega_1}(x) \; {\mbox{ with
}}m(\omega_1)=\alpha\cdot m(\Omega).
$$
It then follows from the unique continuation for heat equation
(\cite{AEWZ})   that $\psi(x,t)=0$. Therefore,
$N\hspace{-1mm}F(\psi)= 0$ if and only if $\psi(x,t)=0$. Finally, a
direct computation shows that
$$
N\hspace{-1mm}F(c\psi)=|c|N\hspace{-1mm}F(\psi), \forall\; \psi\in
Y, \;  c\in \mathbb{R}.
$$
By
$$
\|\beta(\psi_1+\psi_2)\|_{L^2(\Omega\times(s,T))}\le
\|\beta\psi_1\|_{L^2(\Omega\times(s,T))}+\|\beta\psi_2\|_{L^2(\Omega\times(s,T))},\quad \forall\; \beta\in \cB_{s,T},
$$
we have
$$\ba{l}
\disp
\left(\ds\int^T_s \int_\Omega b(x,t)(\psi_1(x,t)+\psi_2(x,t))^2\mathrm
dx \mathrm dt\right)^{1\over 2}\crr\disp \le
\left(\ds\int^T_s \int_\Omega b(x,t)\psi_1(x,t)^2\mathrm
dx \mathrm dt\right)^{1\over 2}+
\left(\ds\int^T_s \int_\Omega b(x,t)\psi_2(x,t)^2\mathrm
dx \mathrm dt\right)^{1\over 2}. \ea$$ So,
$$N\hspace{-1mm}F(\psi_1+\psi_2)\le N\hspace{-1mm}F(\psi_1)+N\hspace{-1mm}F(\psi_2).
$$
This shows that $N\hspace{-1mm}F$ is a norm of  the space $Y$.\endpf

\vspace{3mm}

\bde\label{11.22.1}
 Owing to  Lemma \ref{Le.guo1}, we can denote  the norm given by the
functional $N\hspace{-1mm}F(\cdot)$ as
$\|\cdot\|_{N\hspace{-1mm}F}$. It is clear that the space
$\left(Y,~\|\cdot\|_{N\hspace{-1mm}F}\right)$ is a normed linear
space. We set
$\left(\overline{Y},~\|\cdot\|_{N\hspace{-1mm}F}\right)$ as  the
completion  space of $\left(Y,~\|\cdot\|_{N\hspace{-1mm}F}\right)$.
\ede

Along the same line  in the  proof of Lemma  \ref{lemma2.2}, we have
Lemma \ref{Le.guo2}. \bl\label{Le.guo2} Under an isometric
isomorphism, any element of $ \overline{Y}$ can be expressed as a
function $\hat \f\in C([0,T); L^2(\Omega))$ which satisfies (in the
sense of weak solution)
\begin{equation*}
 \left\{\ba{ll}
  \hat{\varphi}_t(x, t)+\Delta \hat{\varphi}(x,t)-a(x,t)\hat{\varphi}(x,t)=0~~&{\mbox{in }}\Omega\times(s,T),\\
  \hat{\varphi}(x,t)=0 &{\mbox{on }}\partial\Omega\times(s,T),
  \ea\right.
\end{equation*}
and
$N\hspace{-1mm}F(\hat\f)=\lim\limits_{n\rightarrow\infty}N\hspace{-1mm}F(\varphi(\cdot;
z_n))$ for some sequence  $\{z_n\}\subset L^2(\Omega)$, where
$\varphi(\cdot; z_n)$ is the solution of \dref{stateRA} with initial
value $z=z_n$.
\el

We present a further characterization of $\overline{Y}$. \bl~ Let
$Z$ be defined as (\ref{1-1-12}). Then   \be\label{1-3-12} \overline
Y=\left\{\f(\cd;z)\bigm| z\in Z\right\}, \ee where $\f(\cd, z)$ is
the solution to (\ref{stateRA}). Moreover,
\begin{equation}\label{gxy3.41}
\sup\limits_{\psi\in Y} F(b,\psi)=\sup\limits_{\psi\in \overline{Y}}F(b,\psi)
=\sup\limits_{\psi\in \overline{Y_{\beta}}}F(b,\psi).
\end{equation}
\el {\it Proof.} We claim by virtue of  Lemma \ref{Le.guo2} that
\begin{equation}\label{gxy3.40}
\overline{Y}\subseteq L^2(\Omega\times(s,T)).
\end{equation}
Indeed, suppose that $n_0\in\mathbb{N}$ so that $n_0\ge 1/\alpha$.
There are $n_0$ measurable subsets
$\omega_1,\omega_2,\dots,\omega_{n_0}$ of $\Omega$ such that
$$
m(\omega_j)=\alpha\cdot m(\Omega),\quad\forall ~j\in\{1,2,\dots,n_0\},~~\bigcup\limits_{j=1}^{n_0}\omega_j=\Omega.
$$
The inclusion  (\ref{gxy3.40}) then  follows from
\begin{equation}
\label{amss1} \ba{l}
 \ds\int^T_s\int_\Omega\psi(x,t)^2\mathrm dx\mathrm dt
\le\ds\int^T_s\left(\sum\limits_{j=1}^{n_0}\int_\Omega\chi_{\omega_j}(x)\psi(x,t)^2\mathrm dx\right)\mathrm dt \\
\le\ds \sum\limits_{j=1}^{n_0} \ds\int^T_s\int_\Omega\chi_{\omega_j}(x)\psi(x,t)^2\mathrm dx\mathrm dt \le
n_0\|\psi\|_{N\hspace{-1mm}F}^2.
\ea
\end{equation}
Since $\psi(x,t)$ is a generalized function defined on
$\Omega\times(s,T)$ and belongs to $L^2(\Omega\times(s,T))$, and
$\Omega\times{T}$ is the boundary of $\Omega\times(s,T)$, the
inclusion  (\ref{gxy3.40}), together with the trace theorem, implies
(\ref{1-3-12}).  Furthermore, for any $\beta\in\cB_{s,T}$,  by
$$
\|\beta\psi\|_{L^2(s,T: L^2(\Omega))}\le N\hspace{-1mm}F(\psi), \;
\forall\; \psi\in Y,
$$
it follows that \begin{equation}\label{1-3-13}
\overline{Y}\subseteq\overline{Y_{\beta}},\quad
\forall~\beta\in\cB_{s,T}.
\end{equation}
Since  $Y$ is dense in $\overline{Y_{\beta}}$ and
$\sup\limits_{\psi\in Y} F(b,\psi)=\sup\limits_{\psi\in
\overline{Y_{\beta}}}F(b,\psi)$ with $ b=\beta^2$, we obtain
 (\ref{gxy3.41}).
\endpf

Now, we discuss the following  game problem  (with the extend domain of Problem ({RSP2}) or Problem (\ref{gxy3.38})).
\begin{equation}
\begin{array}{ll}
\disp  {\bf Problem (RSP2'):} &\disp \inf\limits_{\psi\in
\overline{Y}}\sup\limits_{
b\in\cB^2_{s,T}}\left[\ds\frac{1}{2}\ds\int^T_s\int_\Omega
b(x,t)\psi(x,t)^2\mathrm dx\mathrm dt+\lan \xi, \psi(s)\ran \right]
\crr & =\disp \inf\limits_{\psi\in
\overline{Y}}\left[\ds\frac{1}{2}\|\psi\|_{N\hspace{-1mm}F}^2+\lan
\xi, \psi(s)\ran \right].
\end{array}
\end{equation}
 Notice that the  functional in Problem (RSP2')  is
strictly convex, coercive, and continuous. Besides, $\overline{Y}$, as a closed subspace
of $L^2(\Omega\times(s, T))$, is also reflexive. Similarly to Lemma
\ref{lemma2.4}, we have  Lemma \ref{lemma3.16}.

 \bl\label{lemma3.16} For any $s\in[0,T)$ and $\xi\in L^2(\Omega)\setminus\{0\}$, Problem  (RSP2') admits
a unique nonzero solution. \el

Now we present the Nash equilibrium problem of two-person zero-sum
game:
\begin{equation}
\begin{array}{ll}
\disp
  {\bf Problem (NEGP):}\;\;  {\mbox{ To find }}\bar b\in\cB^2_{s,T},\,
\bar\psi\in \overline{Y}{\mbox{ such that }}
 F(\bar  b,\bar\psi)&= \sup\limits_{\psi\in \overline{Y}
}F(\bar b,\psi)\crr &=\disp \inf\limits_{ b\in\cB^2_{s,T}}F(
b,\bar\psi),
\end{array}
\end{equation}
where $F(b,\psi)$ is defined by \dref{Taiwan8}. The following
Theorem \ref{theorem3.19} is about  existence  of  Nash
equilibrium to  the  two-person zero-sum game  Problem (NEGP) .

\bt\label{theorem3.19} Let  $\bar\psi$ be a
solution to
 Problem  (RSP2').  Then  Problem (NEGP) admits at least one Nash equilibrium.
  Furthermore, if $\bar\beta$ is a relaxed  optimal actuator location to
   Problem (RP)$^{s,\xi}$, then $(\bar b=\bar\beta^2,
\bar\psi)$ is a Nash equilibrium to Problem (NEGP). Conversely, if
$(\hat  b, \hat\psi)$ is a Nash equilibrium of Problem (NEGP), then
$\hat\psi=\bar\psi$,  and $\hat\beta=\hat b^{1/2}$ is a relaxed
optimal actuator location to  Problem ({RP})$^{s,\xi}$. \et

 \noindent {\bf Proof.}  In terms
of  (\ref{gxy3.41}),
 \be\label{3-15}
 V^+=\inf\limits_{ b\in\cB^2_{s,T}}\sup\limits_{\psi\in Y}F(b,\psi)=\inf\limits_{ b\in\cB^2_{s,T}}\sup\limits_{\psi\in \overline{Y} }F(b,\psi).
\ee
Notice that
$$
V^-=\sup\limits_{\psi\in Y}\inf\limits_{ b\in\cB^2_{s,T}}F(b,\psi)\le \sup\limits_{\psi\in\overline{Y}}\inf\limits_{ b\in\cB^2_{s,T}}F(b,\psi)\le\inf\limits_{ b\in\cB^2_{s,T}}\sup\limits_{\psi\in \overline{Y} }F(b,\psi) .
$$
It follows from Theorem \ref{theorem3.12} that
\begin{equation}\label{gxy3.42}
\inf\limits_{ b\in\cB^2_{s,T}}\sup\limits_{\psi\in \overline{Y} }F(b,\psi)
=\sup\limits_{\psi\in\overline{Y}}\inf\limits_{ b\in\cB^2_{s,T}}F(b,\psi).
\end{equation}
Furthermore, by (\ref{3-15}) and the relation between $\cB_{s,T}$
and $\cB_{s,T}^2$,
 \be\label{3-16}
 \begin{array}{l}
\hbox{if } \bar\beta \hbox{ is a solution to
Problem (RSP)$^{s,\xi}$  },
\hbox{ then }  \bar b \hbox{ is a solution to }
\inf\limits_{ b\in\cB^2_{s,T}}\sup\limits_{\psi\in \overline{Y}
}F(b,\psi);\\
\hbox{if }\bar b \hbox{ is a solution to}
\inf\limits_{ b\in\cB^2_{s,T}}\sup\limits_{\psi\in \overline{Y}
}F(b,\psi),
\hbox{ then }  \bar\beta \hbox{ is a solution to Problem (RSP)$^{s,\xi}$  },
\end{array}
 \ee
 where $\bar b=\bar\beta^2$.
 Recalling  Proposition \ref{prpoperty3.18}, we have
  the following results:
 \begin{itemize}
 \item Equation (\ref{gxy3.42}) ensures that Problem (NEGP) attains its   value;
 \item Problem  (RSP2') admits a unique solution $\bar\psi$ by Lemma \ref{lemma3.16};
 \item Problem (RSP1)   admits  a solution by Theorem  \ref{theoremM1}  and (\ref{3-16}).
  \end{itemize}
  It follows from Proposition
\ref{prpoperty3.18} that Problem (NEGP) admits at least one Nash
equilibrium.
  Furthermore, if $\bar b$ is a solution to $
\inf\limits_{ b\in\cB^2_{s,T}}\sup\limits_{\psi\in \overline{Y}
}F(b,\psi),$
 then $(\bar b,
\bar\psi)$ is a Nash equilibrium to Problem (NEGP). Conversely, if
$(\hat  b, \hat\psi)$ is a Nash equilibrium of Problem (NEGP),
then $\hat b$ is a solution to problem $
\inf\limits_{ b\in\cB^2_{s,T}}\sup\limits_{\psi\in \overline{Y}
}F(b,\psi)$
 and  $ \hat\psi$ solves $\sup\limits_{\psi\in \overline{Y}}\inf\limits_{ b\in\cB^2_{s,T}}
F(b,\psi).$ By the uniqueness from Lemma \ref{lemma3.16}, it holds
that $\hat\psi=\bar\psi$. This, together with (\ref{3-16}) and the
equivalence between Problem ({RSP})$^{s,\xi}$ and Problem
({RP})$^{s,\xi}$,  proves   Theorem \ref{theorem3.19} directly.
\endpf

\section{Proof of the main results}

In this section, we  present proofs for Theorems \ref{theoremM} and
 \ref{theoremM12}.

\subsection{Existence and uniqueness of relaxed optimal actuator location}

Though we have derived the existence for the relaxation problem
(RP)$^{s,\xi}$, existence for  the optimal actuator location to
the classical problem (CP)$^{s,\xi}$ is still not known. To this
purpose, we need to learn more about the  optimal relaxed actuator
location $\bar\beta$. Recall Theorem  \ref{theorem3.19} that  if
$\bar\beta$ is a relaxed actuator location, then $\bar
b=\bar\beta^2$ solves Problem $\sup\limits_{\psi\in \overline{Y}
}F(\bar b,\psi)$. That is to say, $\bar b$ solves \be\label{gxy3.44}
\sup\limits_{ b\in\cB^2_{s,T}}\ds\int^T_s \int_\Omega
b(x,t)\bar\psi(x,t)^2\mathrm dx\mathrm dt. \ee Further, if we denote
$$
 \Gamma=\left\{\gamma\in L^\infty(\Omega;[0,1])\bigm|\int_\Omega\gamma(x)\mathrm dx=\alpha\cdot m(\Omega)\right\},
 $$
then  \be\label{1-4-2} \ds \int_\Omega \bar
b(x,t)\bar\psi(x,t)^2\mathrm dx=\sup\limits_{ \gamma\in\Gamma}\ds
\int_\Omega \gamma(x)\bar\psi(x,t)^2\mathrm dx, \;
 t\in(s,T) \hbox{ a.e.} \ee
 and
 \be\label{4-4-2}  \bar
b(\cd,t)\in \argmax\limits_{ \gamma\in\Gamma}\int_\Omega \gamma(x)\bar\psi(x,t)^2\mathrm dx, \;
 t\in(s,T) \hbox{ a.e.} \ee
 Therefore, we need to discuss
 the following
problem
\be\label{4-38}
\sup\limits_{ \gamma\in\Gamma}\displaystyle\int_\Omega
 \gamma(x)\phi(x)\mathrm dx,
 \ee
  where $\phi\in L^1(\Omega)$.
Similar problem has been discussed in \cite{p2} where $\Gamma$ is replaced by $\cW$. But for the sake of completeness, 
we present here a short proof. 

 Let us define, for any $\phi\in L^1(\Omega)$ and
$c\in\mathbb{R}$,  that
\be\label{01} \ba{l}
\ns\{\phi\ge c\}=\left\{ x\in\Omega\bigm| \phi(x)\ge c\right\},\q\{\phi= c\}=\left\{ x\in\Omega\bigm| \phi(x)= c\right\},\\
\ns\{\phi>c\}=\left\{ x\in\Omega\bigm| \phi(x)>
c\right\},\q\{\phi<c\}=\left\{ x\in\Omega\bigm|
\phi(x)<c\right\}.\ea\ee
 Let \be\label{02}
M_\phi(c)=m(\{\phi\ge c\}) {\mbox{ for any }}\phi\in L^1(\Omega)
\hbox{ and } c\in\mathbb{R}. \ee It is clear that the function
$M_\phi(c)$ is monotone decreasing with respect to $c$. By
$$
\lim\limits_{\ve\rightarrow0+}\{\phi\ge c-\ve\}=\bigcap\limits_{\ve>0}\{\phi\ge c-\ve\}= \{\phi\ge c\},
$$
we have  \be\label{03}
\lim\limits_{\ve\rightarrow0+}M_\phi(c-\ve)=M_\phi(c),\;  \forall\;
c\in\mathbb{R}. \ee This shows that $M_\phi(c)$ is continuous from
the left for any given $\phi\in L^1(\Omega)$. Since
$$
\lim\limits_{c\rightarrow+\infty}M_\phi(c)=0,\q
\lim\limits_{c\rightarrow-\infty}M_\phi(c)=m(\Omega),
$$
the real $c_\phi$ given by
\be\label{04}
 c_\phi=\max\left\{c\in\mathbb{R}\bigm|M_\phi(c)\ge \a\cdot m(\Omega)\right\},
\ee is well-defined. Hence
 \be\label{05} M_\phi(c_\phi)\ge \a\cdot
m(\Omega)\ge M_\phi(c_\phi+)\deq
\lim\limits_{\ve\rightarrow0+}M_{\phi}(c_{\phi}+\ve), \ee and
\be\label{06} M_\phi(c_\phi+\ve)<\alpha m(\Omega),\; \forall\;
\ve>0. \ee
Let \be\label{07}
\bar\a_\phi\deq\displaystyle\frac{M_\phi(c_\phi)}{m(\Omega)},\q
\underline{\a}_\phi\deq\displaystyle\frac{M_\phi(c_\phi+)}{m(\Omega)}.
\ee It follows from (\ref{05}) that \be\label{08} \bar\a_\phi\ge
\a\ge\underline{\a}_\phi. \ee Since
$$
\lim\limits_{\ve\rightarrow0+}\{\phi\ge
c+\ve\}=\bigcup\limits_{\ve>0}\{\phi\ge c+\ve\}= \{\phi> c\},
$$
it follows that
\begin{equation*}
M_{\phi}(c_{\phi}+)=m(\{{\phi}> c_{\phi}\}).
\end{equation*}
By the definition of $\underline{\alpha}_{\phi}$ in \dref{07},
\begin{equation}\label{4-47}
m(\{{\phi}> c_{\phi}\})=\underline{\a}_{\phi}\cdot{m(\Omega)}.
\end{equation}
This,  together with (\ref{07}) and (\ref{08}), implies that
\begin{equation}\label{4-48}
m(\{{\phi}= c_{\phi}\})=(\bar\a_\phi-\underline{\a}_{\phi}){m(\Omega)}\ge (\alpha-\underline{\a}_{\phi}){m(\Omega)}.
\end{equation}
The following result is about problem (\ref{4-38}).
\bl\label{lemma1-4-1}
Let
  $\phi\in W^{1,1}_0(\Omega)$. If $\phi(x)\neq0$ is analytic  in $\Omega$, then Problem (\ref{4-38}) admits a unique solution $\bar\gamma$. Furthermore, it holds that
  \be
  \bar\gamma\in \cW.
  \ee
\el {\it Proof.} Because $\phi(x)$ is analytic, it is clear that
\be\label{1-4-14} m(\{{\phi}= c\})=0 {\mbox{ or }}m(\{{\phi}=
c\})=m(\Omega) \; {\mbox{ for any }}c\in\mathbb{R}. \ee Furthermore,
we claim that \be\label{1-4-15} m(\{{\phi}= c\})=0 \quad{\mbox{for
any }}c\in\mathbb{R}. \ee Indeed, it follows from $\phi\neq 0$ and
(\ref{1-4-14}) that
$$
m(\{{\phi}= 0\})=0.
$$
On the other hand, suppose there is $c\neq 0$ such that
$$
m(\{{\phi}= c\})=m(\Omega).
$$
That is to say, $\phi(x)=c$ in $\Omega$ almost everywhere. Then the
trace of $\phi(x)$ is just  $c$. This contradicts $\phi\in
W^{1,1}_0(\Omega)$. The claim is then  proved.

Let
$c_\phi$, $\bar\a_\phi$, $\underline{\a}_\phi$ defined in (\ref{04}) and (\ref{07}).
It follows from (\ref{4-48}) and (\ref{1-4-15}) that
$$
\bar\a_\phi=\underline{\a}_\phi=\a.
$$

It  follows from (\ref{4-47}) and (\ref{4-48})
 that
 $$m(\{\phi\ge c_\phi\})=\alpha\cdot m(\Omega).$$
 That is, $\{\phi\ge c_\phi\}\in\cW.$
Since  $\Gamma$ is  the convex hull of
$\{\chi_\omega|\omega\in\cW\}$, it holds
$$
\sup\limits_{\gamma\in \Gamma}\int\gamma \phi\mathrm
dx=\sup\limits_{\omega\in \cW}\int\chi_\omega \phi\mathrm dx.
$$
If we can show that
$$
\int\chi_{\{\phi\ge c_\phi\}} \phi\mathrm dx>\int
\chi_\omega\phi\mathrm dx, \; \forall\;
 \omega\in\cW,\,\chi_\omega\neq\chi_{\{\phi\ge c_\phi\}},
$$
then  $\chi_{\{\phi\ge c_\phi\}}$ is the unique solution to problem
(\ref{4-38}) and  belongs to $\cW$. To this purpose,  let
$\omega_1=\omega\setminus\{\phi\ge c_\phi\}$,  $\omega_2=\{\phi\ge
c_\phi\}\setminus\omega$, and $\omega_3=\omega\cap\{\phi\ge
c_\phi\}$. Since  $\omega$ and $\{\phi\ge c_\phi\}$ belong to $\cW$,
it holds
$$
m(\omega_1)=m(\omega_2)\neq 0.
$$
On the other hand, since
$$
\phi(x)\ge c_\phi>\phi(y)\; \forall\; x\in\omega_2,\, y\in\omega_1,
$$
we thus have
$$
\int\chi_{\{\phi\ge c_\phi\}} \phi\mathrm
dx=\int_{\omega_2}\phi\mathrm dx+\int_{\omega_3}\phi\mathrm
dx>\int_{\omega_1}\phi\mathrm dx+\int_{\omega_3}\phi\mathrm dx=\int
\chi_\omega\phi\mathrm dx.
$$
Therefore,  $\chi_{\{\phi\ge c_\phi\}}$ is the unique solution to
problem (\ref{4-38}) and  belongs to $\cW$.
\endpf

\bigskip

\noindent {\bf  Proof of Theorem \ref{theoremM}}.  Recall that the
coefficient $a(x,t)$ is analytic. Thus  the solution to  Equation
(\ref{stateRA}) with the initial condition $z\in L^2(\Omega)$ is
also analytic in $\Omega\times(s,T)$ (\cite{AEWZ}).  As  the
solution to  Problem (RSP2'),
$$\bar\psi(\cd,T-\varepsilon)\in L^2(\Omega) \; {\mbox{ for any }}\varepsilon>0.
$$
Thus $\bar\psi$ is analytic in $\Omega\times(s,T-\varepsilon)$. By
the arbitrariness of $\varepsilon$,  $\bar\psi$ is analytic in
$\Omega\times(s,T)$. On the other hand,  it follows from the smooth
effect of the heat equation that
$$
\bar\psi(\cd,t)\in H^1_0(\Omega)\; \mbox{ for any }t\in(s,T).
$$
Those, together with the non-singularity of $\bar\psi$, imply that
$$
\bar\psi(\cd,t)^2{\mbox{ is nonzero analytic in }} \Omega {\mbox{ and }} \bar\psi(\cd,t)^2\in W^{1,1}_0(\Omega) {\mbox{ for any }}t\in(s,T).
$$
By  Lemma \ref{lemma1-4-1} and (\ref{1-4-2}),  $\bar b$ is unique
and  belongs to $\cW_{s,T}$. Therefore, it follows from Theorem
\ref{theorem3.19} that any relaxed optimal actuator location must be
classical and unique. We thus complete the proof of the theorem.
\endpf

\vspace{0.3cm}

\noindent {\bf Proof of Theorem \ref{theoremM12}}. We  use the
synthetic method, to obtain  the feedback and prove the
corresponding result by the dynamic programming approach. The
synthetic method is a method to be used to construct a feedback
control through open-loop control reflected mathematically by
\dref{G4.17} and \dref{G4.18} later (see, e.g., \cite{Wnew}).

Now, for any $(s,\xi)\in[0,T)\times L^2(\Omega)\setminus\{0\}$, denote the optimal
actuator location by $w^{s,\xi}\in\cW_{s,T}$ and the corresponding
optimal control of Problem (NP)$^{s,\xi}_{w^{s,\xi}}$ by $u^{s,\xi}\in L^2(\Omega\times(s,T))$.
Write the
corresponding optimal trajectory by $y^{s,\xi}\in
C([s,T];L^2(\Omega))$. Based on these notations, we begin to define $\cF:[0,T)\times L^2(\Omega)\mapsto
\cW$ by
\begin{equation}\label{G4.17}  \cF(s,\xi)=w^{s,\xi}(s)
\; \mbox{ for any }(s,\xi)\in[0,T)\times L^2(\Omega),
\end{equation}
 and  define $\cG: [0,T)\times L^2(\Omega)\mapsto L^2(\Omega)$ by
 \begin{equation}\label{G4.18}
\cG(s,\xi)=u^{s,\xi}(s) \; \mbox{ for any }(s,\xi)\in[0,T)\times
L^2(\Omega).
\end{equation}
The above definition is well-defined. Indeed, as the solution of Problem (RSP2'), $\bar\psi\in C([s, T);L^2(\Omega))$. It follows from (\ref{4-4-2}) that
$w^{s,\xi}\in C([s, T);L^2(\Omega))$. By Lemma \ref{lemma2.4},
$u^{s,\xi}\in w^{s,\xi}\overline Y_{w^{s,\xi}}$. This, together with the continuity of $w^{s,\xi}$, implies  that
$u^{s,\xi}\in C([s, T);L^2(\Omega))$.

\medskip

Fix $(s,\xi)\in[0,T)\times L^2(\Omega)$. We will show that  $y^{s,\xi}$ defined as above is just the unique solution of Equation (\ref{1-1-4})  satisfying
$y^{\cF,\cG}((x,T);s, \xi)=0$ and (\ref{1-1-5})-(\ref{1-1-7}). The proof will be carried out by the following several steps.

\medskip

{\it  Step 1:  $u^{s,\xi}\bigl|_{[t,T)}$ is the solution to Problem (NP)$^{t, y^{s,\xi}(t)}_{w^{s,\xi}\bigm|_{[t,T)}}$.}

Notice that
$$y\left(T;\,w^{s,\xi}\bigm|_{[t,T)}, u^{s,\xi}\bigl|_{[t,T)}; t, y^{s,\xi}(t)\bigl|_{[t,T)}\right)=
y\left(T;\,w^{s,\xi}, u^{s,\xi};s,\xi\right)=0.
$$
If there is $v\in L^2(\Omega\times(t,T)$ such that $$
\ba{l}
y\left(T;\,w^{s,\xi}\bigm|_{[t,T)}, v; t, y^{s,\xi}(t)\bigl|_{[t,T)}\right)=0
{\mbox{ with }}
\|v\|_{ L^2(\Omega\times(t,T))}<\left\|u^{s,\xi}\bigl|_{[t,T)}\right\|_{L^2(\Omega\times(t,T))},
\ea
$$
by setting
$$
\ba{l}
\hat v(r)=\left\{\ba{ll}
u^{s,\xi}(r),\hspace{5pt} &{\mbox{when }}r\in[s,t)\\
v(r)&{\mbox{when }}r\in[t,T),
\ea\right.
\ea
$$
we  find that $\hat v\in L^2(\Omega\times(s,T))$ satisfies  $y\left(T;\,w^{s,\xi}, \hat v;s,\xi\right)=0$
and $\|\hat v\|_{ L^2(\Omega\times[s,T))}<\left\|u^{s,\xi}\right\|_{L^2(\Omega\times[s,T))}$.
This means that $\hat v$ solves Problem (NP)$^{t, y^{s,\xi}(t)}_{w^{s,\xi}}$, which contradicts
with the optimality  of $u^{s,\xi}$ and thus leads to claim of step 1.

\medskip

{\it Step 2:  $w^{s,\xi}\bigl|_{[t,T)}$ is the solution to Problem (CP)$^{t, y^{s,\xi}(t)}$.}

Assume the above claim is false. Then there is   $\hat  w\in \cW_{t,T}$ solving Problem (CP)$^{t, y^{s,\xi}(t)}$. Denote by $\tilde v\in L^2(\Omega\times(t,T))$ the solution to Problem (NP)$^{t, y^{s,\xi}(t)}_{\hat w}$.
By setting
$$
\tilde  w(r)=\left\{\ba{ll}
w^{s,\xi}(r),\hspace{5pt} &{\mbox{when }}r\in[s,t)\\
\hat w(r)&{\mbox{when }}r\in[t,T),\ea
\right.
\quad
 \hat v(r)=\left\{\ba{ll}
u^{s,\xi}(r),\hspace{5pt} &{\mbox{when }}r\in[s,t)\\
\tilde v(r)&{\mbox{when }}r\in[t,T),\ea
\right.
$$
we   find that $y\left(T;\,\tilde w, \hat v;s,\xi\right)=0$. Now we claim
\be\label{4-4-21}
\|\hat v\|_{ L^2(\Omega\times[s,T))}<\left\|u^{s,\xi}\right\|_{L^2(\Omega\times[s,T))}.
\ee
Indeed, by {\it Step 1},
$$N\left(w^{s,\xi}\big|_{[t,T)}; t, y^{s,\xi}(t)\right)=
\left\|u^{s,\xi}\big|_{[t,T)}\right\|_{L^2(\Omega\times(t,T))}.
$$
Because
$N\left(\hat w; t, y^{s,\xi}(t)\right)=
\left\|\tilde v\right\|_{L^2(\Omega\times(t,T))}
$, it follows from the unique optimality of $\hat w$ that
$$
\left\|u^{s,\xi}\big|_{[t,T)}\right\|_{L^2(\Omega\times(t,T))}
>\left\|\tilde v\right\|_{L^2(\Omega\times(t,T))}.
$$
This implies (\ref{4-4-21}) and hence $\tilde w$ solves   Problem (CP)$^{s,\xi}$ which is impossible. We
Thus conclude the claim of {\it Step 2}.

\medskip

{\it Step 3: $y^{s,\xi}$ is the unique solution to (\ref{1-1-4}) satisfying
 $y^{s,\xi}(T)=0$ and (\ref{1-1-5})-(\ref{1-1-7}).}

 It is clear that   $y^{s,\xi}(T)=0$. From {\it Step 2}, we have
\be\label{4-4-22}
\cF(t, y^{s,\xi}(t))=w^{s,\xi}(t)\quad{\mbox{for any }}t\in[s,T).
\ee
By {\it Step 1}, we have
\be\label{4-4-23}
\cG(t, y^{s,\xi}(t))=u^{s,\xi}(t)\quad{\mbox{for any }}t\in[s,T).
\ee
Thus
 $y^{s,\xi}$ is a solution to (\ref{1-1-4}).  In addition,
 it follows from (\ref{4-4-22})-(\ref{4-4-23}) and the definition of  $w^{\cF,\cG}(s,\xi)$ and $u^{\cF,\cG}(s,\xi)$ that
\be\label{4-4-24}
 w^{\cF,\cG}(s,\xi)=w^{s,\xi}\in \cW_{s,T},\quad u^{\cF,\cG}(s,\xi)=u^{s,\xi}\in L^2(\Omega\times(s,T).
 \ee
This gives (\ref{1-1-6})-(\ref{1-1-7}). The identities  (\ref{1-1-5}) follow straightforwardly  from the optimality of $u^{s,\xi}$.
Therefore, $y^{s,\xi}$ is a solution to (\ref{1-1-4}) satisfying
 (\ref{1-1-5})-(\ref{1-1-7}).

Finally, we come up uniqueness. From (\ref{1-1-4}), we find  that
 $w^{\cF,\cG}(s,\xi)$ is a solution to Problem (CP)$^{s,\xi}$.  The identities (\ref{4-4-24}) follow from the uniqueness.
In addition,  as the solution to  (\ref{1-1-4}),$y^{s,\xi}$  must satisfy the following equation
$$
 \left\{\ba{ll}
  y_t(x,t)-\Delta y(x,t)+a(x,t)y(x,t)=\left(w^{s,\xi}u^{s,\xi} \right)(x,t)~~&{\mbox{in }}\Omega\times(s,T),\\
 y(x,t)=0 &{\mbox{on }}\partial\Omega\times(s,T),\\
  y(x,s)=\xi(x)&{\mbox{in }}\Omega,
  \ea\right.
$$
It is clear that $y^{s,\xi}$ is the unique solution.
We thus complete the proof of the theorem.\endpf

\end{document}